\documentclass[11pt,a4paper,reqno]{amsart}
\usepackage{stmaryrd}
\usepackage[english]{babel}
\usepackage{pst-grad} 
\usepackage{pst-plot} 
\usepackage{pstricks}
\usepackage{amsmath,amssymb,mathrsfs,amsthm, mathtools}
\usepackage{tikz} 
\usepackage{mathcomp,wasysym}  
\usepackage{graphicx}  
\usepackage[all]{xy} \xyoption{arc} \xyoption{color}
\usepackage{epsfig}
\usepackage{cite}

\newcommand{\pare}[1]{\left( #1 \right)}

\newcommand{\cH}{\mathcal{H}}

\normalsize
\normalsize
\def\comm#1#2{{\left\llbracket#1,#2\right\rrbracket}}

\usepackage[bookmarks,colorlinks]{hyperref}

\newtheorem{theorem}{Theorem}

\theoremstyle{definition}
\newtheorem{remark}{Remark}

\title{Well-posedness of a water wave model with viscous effects}

\author[R. Granero-Belinch\'{o}n]{Rafael Granero-Belinch\'{o}n}
\address{Departamento  de  Matem\'aticas,  Estad\'istica  y  Computaci\'on,  Universidad  de Cantabria.  Avda.  Los  Castros  s/n,  Santander,  Spain.}
\email{rafael.granero@unican.es}

\author[S. Scrobogna]{Stefano Scrobogna}
\address{IMUS, Av. de la Reina Mercedes, s/n, 41012 Sevilla, Spain \& BCAM, Mazarredo 14, 48009, Bilbao, Spain}
\email{sscrobogna@bcamath.org}

\begin{document}
\begin{abstract}
Starting from the paper by Dias, Dyachenko and Zakharov (\emph{Physics Letters A, 2008}) on viscous water waves, we derive a model that describes water waves with viscosity moving in deep water with or without surface tension effects. This equation takes the form of a nonlocal fourth order wave equation and retains the main contributions to the dynamics of the free surface. Then, we prove the well-posedness in Sobolev spaces of such equation.
\end{abstract}



\maketitle
{\small
\tableofcontents}

\allowdisplaybreaks
\section{Introduction}
The motion of the free surface of an incompressible fluid has been studied for centuries \cite{Stokes_1847}. Many of these works examine the case where the fluid is inviscid and irrotational and thus they have to study the Euler equations. However, there are situations where the viscous damping cannot be neglected and should be considered. The first works in this direction date back to Boussinesq \cite{boussinesq1895lois} and Lamb \cite{lamb1932hydrodynamics}. Then, Ruvinsky \& Freidman \cite{ruvinsky1987fine} formulated a system of equations for weakly damped water waves in deep water (see also \cite{ruvinsky1985improvement,ruvinsky1991numerical,khariff1996frequency}). Similar models were also proposed and studied by Longuet-Higgins \cite{longuet1992theory}, Jiang, Ting, Perlin \& Schultz \cite{jiang1996moderate}, Joseph \& Wang \cite[Equation (6.7) and (6.8)]{joseph2004dissipation}, Wang \& Joseph \cite{wang2006purely} and Wu, Liu \& Yue \cite{wu2006note}. The interested reader can refer to \cite{graneroscrobo} for more details on these contributions. 

More recently, Dias, Dyachenko \& Zakharov \cite{dias2008theory} proposed the following system
\begin{subequations}\label{eq:DDZ}
\begin{align}
\Delta \phi&=0&&\text{ in }\Omega(t),\\
\rho\left(\phi_t+\frac{1}{2}|\nabla\phi|^2+ Gh\right)&=-2\mu \partial_2^2\phi&&\text{ on }\Gamma(t),\\
h_t&=\nabla\phi\cdot \left(-\partial_{1}h,1\right)+2\frac{\mu}{\rho}\partial_{1}^2h &&\text{ on }\Gamma(t),
\end{align}
\end{subequations}
as a model of viscous water waves. Here the $t\in[0,T]$ is the time while $\Omega(t)$ and $\Gamma(t)$ are the the region occupied by the fluid and the water wave. Furthermore, $h$ denotes the height of the interface and $\phi$ is the velocity potential.

This model obtained a considerable amount of attention and was also considered by several other authors. For instance, Dutykh \& Dias \cite{dutykh2007viscous} deduced a new system taking into account the effects of the bottom topography (see also \cite{dutykh2009visco, dutykh2007dissipative,dutykhA,dutykh2009visco}). Kakleas \& Nicholls \cite{kakleas2010numerical}, using the analytic dependence of the Dirichlet-Neumann operator, derived and asymtotic model of \eqref{eq:DDZ}. This asymptotic model takes the form of a system of two equations and is the viscous analog of the classical Craig-Sulem WW2 model \cite{craig1993numerical}. For the well-posedness of these systems we refer to the works Ambrose, Bona \& Nicholls \cite{ambrose2012well}, Ngom \& Nicholls \cite{ngom2018well} or the new preprint by the authors \cite{graneroscrobo2}.

In the very recent paper \cite{graneroscrobo}, the authors derived two different asymptotic models in the regime of small steepness. In particular, starting from the model for water waves with viscosity proposed by Dias, Dyachenko and Zakharov \eqref{eq:DDZ}, the authors obtained a nonlocal wave equation for the free surface. This wave equation contains terms at different scales. The purpose of this paper is to propose a model that retains only the main contributions and to study whether this new model is well-posed in the sense of Hadamard. In particular, this note is devoted to the derivation and mathematical study of the following nonlocal wave equation
\begin{equation}\label{NEW}
\begin{aligned}
&\begin{multlined}
f_{tt}+2\delta\Lambda^2 f_t+ \Lambda f+\beta\Lambda^3 f+\delta^2\Lambda^4 f= -\Lambda\left(\left(\cH f_t\right)^2\right)+\partial_x\comm{\cH}{f}\Lambda f \\+\beta\partial_x\comm{\cH}{f}\Lambda^3 f
+\delta\partial_x\comm{\cH}{\cH f _t}\cH \partial_x ^2 f+\delta\Lambda\left(\cH f _t\cH \partial_x ^2 f \right)\\-\delta\partial_x\comm{\partial_x^2}{f} \cH f _t , 
\end{multlined},
\end{aligned}
\end{equation}
where $(x,t)\in\mathbb{S}^1\times[0,T]$, $\mathbb{S}^1$ denotes the one-dimensional flat torus (identified with $[-\pi,\pi]$ equipped with periodic boundary conditions), $\delta>0$ is a dimensionless parameter reflecting the viscous damping of water waves, $\beta\geq0$ is the Bond number, the operators $\mathcal{H}$ and $\Lambda$ stand for the Hilbert transform and the Calder\'on operator
\begin{align}\label{Hilbert}
\widehat{\mathcal{H}f}(k)=-i\text{sgn}(k) \hat{f}(k) \,, \ \ 
\widehat{\Lambda f}(k)=|k|\hat{f}(k)\,, \ \ \widehat{\Lambda^s f}(k)&=|k|^s\hat{f}(k)\,,
\end{align}
and 
$$
\comm{A}{B}f=A(Bf)-B(Af),
$$ is the commutator between two operators acting on the function $f$. Together with \eqref{NEW}, we have to consider the initial data
\begin{align}\label{initialdata}
f(x,0)=f_0(x), && f_t(x,0)=f_1(x).
\end{align}
From this point onwards and without losing generality, we will assume that the initial data have zero mean. In fact, this  equation is obtained as a simplified mathematical model of water waves with viscosity (see below for a more detailed explanation). The inviscid analog of equation \eqref{NEW} was studied by different authors: Matsuno \cite{matsuno1992nonlinear,matsuno1993two,matsuno1993nonlinear}, Akers \& Milewski \cite{AkMi2010} and Cheng, Granero-Belinch\'on, Shkoller and Wilkening \cite{aurther2019rigorous} derived the inviscid model. Furthermore, Cheng, Granero-Belinch\'on, Shkoller and Wilkening also proved that it is well-posed in a class of analytic functions and performed numerical simulations. Finally, let us mention that Akers \& Nicholls \cite{AkNi2010} studied the inviscid model under the traveling solitary wave ansatz. Finally, we note that there is a striking similarity between the inviscid version of \eqref{NEW} and models describing water waves in a porous media \cite{granero2019asymptotic,granero2018asymptotic}.

The main result of this work is the following theorem
\begin{theorem}\label{th}Fix $\delta>0$ and $\beta\geq0$. Let $(f_0,f_1)\in H^{5.5}\times H^{3.5}$ be the arbitrary initial data \eqref{initialdata} for \eqref{NEW}. Then, there exist $0<T\leq\infty$ such that there is a unique solution
\begin{align*}
f & \in C\left([0,T],H^{5.5}\right), \\
f_t & \in C\left([0,T],H^{3.5}\right)\cap L^2(0,T;H^{4.5}).
\end{align*}
for \eqref{NEW} for a short enough lifespan $0<T\leq \infty$.
\end{theorem}
The proof of this theorem relies on appropriate energy estimates and a Picard iteration method for a mollified version of \eqref{NEW} (see \cite{granero2017model} for more details). We would like to emphasize that, in order the initial data can be arbitrarily big, the energy estimates need to exploit the fine structure of some terms in the nonlinearity of \eqref{NEW}. Let us also explain the reason behind the use of the space $H^{5.5}$ (once this spaces is fixed, the other ones are also fixed due to scaling): as we are interested in classical solutions, we need at least $H^\alpha$ with $\alpha>4.5$. That means testing against $\Lambda^{\sigma}f_t$ with $\sigma>5$ ($\sigma$ and $\alpha$ are related via $\alpha=(\sigma+4)/2$). For the sake of simplicity, we consider $\sigma=7$ and that leads to the space $H^{5.5}$ in the statement. However, with more  sophisticates estimates we could also reach the space $H^5$.

\section{Motivation}
Let us briefly sketch the derivation of \eqref{NEW}. Starting from the model for water waves with viscosity system proposed by Dias, Dyachenko, and Zakharov \cite{dias2008theory}, we derived the following nonlocal wave equation \cite{graneroscrobo} 
\begin{multline}\label{models2}
f_{tt}-(\alpha_1+\alpha_2)\partial_{1}^2 f_t+ \Lambda f+\beta\Lambda^3 f+\alpha_1\alpha_2\partial_{1}^4 f\\
= \varepsilon\bigg\lbrace-\Lambda\left(\left(\cH f_t\right)^2\right)+\partial_x\comm{\cH}{f}\Lambda f +\beta\partial_x\comm{\cH}{f}\Lambda^3 f
\\
+\alpha_2\partial_x\comm{\cH}{\cH f _t}\cH \partial_x ^2 f+\alpha_2\Lambda\left(\cH f _t\cH \partial_x ^2 f \right)
+\alpha_1\alpha_2\partial_x\comm{\partial_x^2}{f}\Lambda\partial_{1} f\\
-\alpha_1\partial_x\comm{\partial_x^2}{f} \cH f _t
-\alpha_2\alpha_2\partial_x\comm{\cH}{\partial_{1}^2 f}\partial_{1}^2 f \bigg\rbrace. 
\end{multline}

Roughly speaking, some typical values of the dimensionless parameters are (see \cite{graneroscrobo}):
\begin{align} \label{eq:dimensionless_parameters1}
\varepsilon\approx O(10^{-2}), &&  \beta\approx O(10^{-5}),  && \alpha_1,\alpha_2\approx O(10^{-4}).
\end{align}

Thus, we see that the smallest terms are those with $\varepsilon\alpha_2\alpha_2,\varepsilon\alpha_2\alpha_1.$ These terms are $O(10^{-10})$. It is reasonable to think that they are not crucial for the evolution of the interface and that, even neglecting them, we still can describe most of the dynamics. As a consequence, we obtain the new model (where $\alpha_1=\alpha_2=\delta$) 
\begin{multline}\label{NEWB}
f_{tt}+2\delta\Lambda^2 f_t+ \Lambda f+\beta\Lambda^3 f+\delta^2\Lambda^4 f= \varepsilon\bigg\lbrace-\Lambda\left(\left(\cH f_t\right)^2\right)+\partial_x\comm{\cH}{f}\Lambda f \\+\beta\partial_x\comm{\cH}{f}\Lambda^3 f
+\delta\partial_x\comm{\cH}{\cH f _t}\cH \partial_x ^2 f+\delta\Lambda\left(\cH f _t\cH \partial_x ^2 f \right)\\-\delta\partial_x\comm{\partial_x^2}{f} \cH f _t \bigg\rbrace,\;(x,t)\in\mathbb{S}^1\times[0,T], 
\end{multline}
As our analysis is independent of the steepness, without losing generality, we take the steepness parameter to be identically one, \emph{i.e.} $\varepsilon=1$. As a consequence, we recover \eqref{NEW}.

\section{Proof of Theorem \ref{th}}
We are going to find appropriate \emph{a priori} bounds for the energy
\begin{equation}\label{eq:energy}
\mathfrak{E}(t)=\max_{0\leq s\leq t}\left\{\|f(s)\|_{\dot{H}^{4}}^2+\beta \|f(s)\|_{\dot{H}^{5}}^2+\delta^2 \|f(s)\|_{\dot{H}^{5.5}}^2+\|f_t(s)\|_{\dot{H}^{3.5}}^2\right\}.
\end{equation}
It is also convenient to define the dissipation
\begin{equation}\label{eq:dissipation}
\mathfrak{D}(t)=2\delta\|f_t(t)\|_{\dot{H}^{4.5}}^2.
\end{equation}
Our goal is to find a polynomial inequality of the form
$$
\mathfrak{E}(t)\leq \mathcal{M}_0(f_0,f_1)+t^{\alpha}P(\mathfrak{E}(t)),
$$
for certain $\alpha>0$, constant $\mathcal{M}_0(f_0,f_1)$ and polynomial $P$. Testing \eqref{NEW} against $\Lambda^7f_t$ and integrating by parts, we find that
$$
\frac{1}{2}\mathfrak{E}(t)+\int_0^t\mathfrak{D}(s)\ \text{d} s=\frac{1}{2}\mathfrak{E}(0)+\sum_{i=1}^6 I_i,
$$
with 
\begin{subequations}\label{eq:Ij}
\begin{align}
I_1&=-\int_0^t\int_{\mathbb{S}^1}\Lambda\left(\left(\cH f_t\right)^2\right)\Lambda^7 f_t \ \text{d} x \ \text{d}s\\
I_2&=\int_0^t\int_{\mathbb{S}^1}\partial_x\comm{\cH}{f}\Lambda f\Lambda^7 f_t \ \text{d} x \ \text{d}s\\
I_3&=\beta\int_0^t\int_{\mathbb{S}^1}\partial_x\comm{\cH}{f}\Lambda^3 f\Lambda^7 f_t \ \text{d} x \ \text{d}s\\
I_4&=\delta\int_0^t\int_{\mathbb{S}^1}\partial_x\comm{\cH}{\cH f _t}\cH \partial_x ^2 f\Lambda^7 f_t \ \text{d} x \ \text{d}s\\
I_5&=\delta\int_0^t\int_{\mathbb{S}^1}\Lambda\left(\cH f _t\cH \partial_x ^2 f \right)\Lambda^7 f_t \ \text{d} x \ \text{d}s\\
I_6&=-\delta\int_0^t\int_{\mathbb{S}^1}\partial_x\comm{\partial_x^2}{f} \cH f _t\Lambda^7 f_t \ \text{d} x \ \text{d}s.
\end{align}
\end{subequations}
Using the self-adjointness of the operator $\Lambda$ together with H\"older's inequality and the Sobolev embedding
$$
\|g\|_{L^4}\leq C\|g\|_{H^{0.25}},
$$
we find that
\begin{align}
I_1&=-\int_0^t\int_{\mathbb{S}^1}\left(\left(\cH f_t\right)^2\right)\Lambda^8 f_t \ \text{d} x \ \text{d}s\nonumber\\
&=-\int_0^t\int_{\mathbb{S}^1}\left(\left(\cH f_t\right)^2\right)\partial_x^4 \Lambda^4 f_t \ \text{d} x \ \text{d}s\nonumber\\
&=-\int_0^t\int_{\mathbb{S}^1}\partial_x^4\left(\left(\cH f_t\right)^2\right)\Lambda^4 f_t \ \text{d} x \ \text{d}s\nonumber\\
&=-\int_0^t\int_{\mathbb{S}^1}\left(2\cH f_t\Lambda\partial_x^3 f_t+6(\Lambda\partial_x f_t)^2+8\Lambda f_t\partial_x^2\Lambda f_t\right)\Lambda^4 f_t \ \text{d} x \ \text{d}s\nonumber\\
&\leq C\int_0^t\|f_t\|_{\dot{H}^4}\left(\|f_t\|_{\dot{H}^4}\|\cH f_t\|_{L^\infty}+\|f_t\|_{\dot{H}^{2.25}}^2+\|f_t\|_{\dot{H}^3}\|\Lambda f_t\|_{L^\infty}\right)\ \text{d} s\nonumber\\
&\leq C\int_0^t\|f_t\|_{\dot{H}^{4.5}}\|f_t\|_{\dot{H}^{3.5}}\|f_t\|_{\dot{H}^{2.25}}\ \text{d} s.	\label{ineq:I1}
\end{align}
Furthermore, using interpolation between Sobolev spaces, the embedding
$$
H^s\subset H^{r}\,,r\leq s, 
$$
the continuity of the Hilbert transform
\begin{align*}
\|\cH F\|_{L^p}\leq C\|F\|_{L^p},  && 1< p<\infty,
\end{align*}
and Young's inequality, we can further compute that
\begin{align}
I_1&\leq \frac{C}{\delta}t\mathfrak{E}^2+\frac{\delta}{4}\int_0^t\mathfrak{D}\ \text{d} s\label{I12}
\end{align}
We recall the following commutator estimate (see equation (1.13) in \cite{dawson2008decay})
\begin{align}\label{commutatorH}
\|\partial_x^\ell \comm{\cH}{u}\partial_x^m v\|_{L^p}\leq C\|\partial_x^{\ell+m}u\|_{L^\infty}\|v\|_{L^p}, && p\in(1,\infty), && \ell,m\in\mathbb{N}.
\end{align}
Equipped with \eqref{commutatorH}, we can estimate $I_2$ as follows
\begin{align}
I_2&=\int_0^t\int_{\mathbb{S}^1}\Lambda^3\partial_x\comm{\cH}{f}\Lambda f\Lambda^4 f_t \ \text{d} x \ \text{d}s\nonumber\\
&\leq\int_0^t\|\partial_x^4\comm{\cH}{f}\Lambda f\|_{L^2}\|\Lambda^4 f_t\|_{L^2} \ \text{d} s\nonumber\\
&\leq C\int_0^t\|\partial_x^4f\|_{L^\infty}\|\Lambda f\|_{L^2}\|\Lambda^4 f_t\|_{L^2} \ \text{d} s\nonumber\\
&\leq C\int_0^t\|\partial_x^5f\|_{L^2}\|\Lambda f\|_{L^2}\|\Lambda^4 f_t\|_{L^2} \ \text{d} s\nonumber\\
&\leq \frac{C}{\delta^3}t\mathfrak{E}^2\pare{t}+\frac{\delta}{4}\int_0^t\mathfrak{D}\pare{s}\ \text{d} s\label{I22}.
\end{align}
Using \eqref{commutatorH}, we have that $I_3$ can be estimated as follows:
\begin{align}
I_3&=\int_0^t\int_{\mathbb{S}^1}\beta\Lambda^3\partial_x\comm{\cH}{f}\Lambda^3 f\Lambda^{4} f_t \ \text{d} x \ \text{d}s\nonumber\\
&\leq C\beta\int_0^t\|f_t\|_{\dot{H}^{4}}\|\partial_x^4 f\|_{L^\infty}\|f\|_{\dot{H}^3}\ \text{d} s\nonumber\\
&\leq \frac{C}{\delta}t\mathfrak{E}^2+\frac{1}{4}\int_0^t\mathfrak{D}\ \text{d} s\label{I32}.
\end{align}
We can decompose $I_4$ as follows
\begin{align*}
I_4&=\delta\int_0^t\int_{\mathbb{S}^1}\left[\Lambda(\cH f_t\Lambda\partial_x f)+\partial_x(\cH f_t\partial_x^2 f)\right]\Lambda^{7} f_t \ \text{d} x \ \text{d}s\\
&=\delta\int_0^t\int_{\mathbb{S}^1}\left[\Lambda(\cH f_t\Lambda\partial_x f)-\partial_x(\cH f_t\Lambda^2 f)\right]\Lambda^{7} f_t \ \text{d} x \ \text{d}s\\
&=J_1^4+J_2^4,
\end{align*}
with
\begin{align*}
J_1^4&=\delta\int_0^t\int_{\mathbb{S}^1}\Lambda(\cH f_t\Lambda\partial_x f)\Lambda^{7} f_t \ \text{d} x \ \text{d}s\\
J_2^4&=-\delta\int_0^t\int_{\mathbb{S}^1}\partial_x(\cH f_t\Lambda^2 f)\Lambda^{7} f_t \ \text{d} x \ \text{d}s.
\end{align*}
We will use the fractional Leibniz rule (see \cite{grafakos2014kato,kato1988commutator,kenig1993well}):
$$
\|\Lambda^s(uv)\|_{L^p}\leq C\left(\|\Lambda^s u\|_{L^{p_1}}\|v\|_{L^{p_2}}+\|\Lambda^s v\|_{L^{p_3}}\|u\|_{L^{p_4}}\right),
$$
which holds whenever
$$
\frac{1}{p}=\frac{1}{p_1}+\frac{1}{p_2}=\frac{1}{p_3}+\frac{1}{p_4}\qquad \mbox{where $1/2<p<\infty,1<p_i\leq\infty$},
$$
and $s>\max\{0,1/p-1\}$. Recalling the inequality (see for instance \cite{cheng2016well} for a proof)
\begin{equation}\label{almostalgebra}
\|uv\|_{H^{0.5}}\leq C_\sigma\|u\|_{H^{0.5}}\|v\|_{H^{0.5+\sigma}},\;\forall\sigma>0,
\end{equation}
using the duality pairing $H^{0.5}-H^{-0.5}$, the Sobolev embedding, the fractional Leibniz rule and the self-adjointness of the operator $\Lambda$, we compute
\begin{align*}
J_1^4&=\delta\int_0^t\int_{\mathbb{S}^1}\Lambda^3(\cH f_t\Lambda\partial_x f)\Lambda^{5} f_t \ \text{d} x \ \text{d}s\nonumber\\
&\leq\delta C\int_0^t\|\Lambda^{3.5}(\cH f_t\Lambda\partial_x f)\|_{L^2}\|\Lambda^{5} f_t\|_{H^{-0.5}}\ \text{d} s\nonumber\\
&\leq\delta C\int_0^t(\|f_t\|_{\dot{H}^{1}}\|f\|_{\dot{H}^{5.5}}+\|f_t\|_{\dot{H}^{3.5}}\|f\|_{\dot{H}^{3}})\|f_t\|_{H^{4.5}}\ \text{d} s\nonumber\\
&\leq \frac{C}{\delta}t\mathfrak{E}^2+\frac{\delta}{4}\int_0^t\mathfrak{D}\ \text{d} s.
\end{align*}
The terms $J^4_2$ and $I_5=J^4_1$ can be estimated in a similar way and we find that
\begin{equation}\label{I42}
I_4+I_5\leq \frac{C}{\delta}t\mathfrak{E}^2 \pare{t} +\frac{\delta}{2}\int_0^t\mathfrak{D}\pare{s} \ \text{d} s. 
\end{equation}
Now we are left with $I_6$. This is the more delicate term due to the high number of derivatives present in the terms involving $f_t$. In fact, if we just use the $H^{0.5}-H^{-0.5}$ duality, we would get an estimate of the form
$$
I^6\leq 2\delta\int_0^t\|f\|_{\dot{H}^{2}}\|f_t\|_{\dot{H}^{4.5}}^2\ \text{d} s+\text{l.o.t.}\leq 2\delta \sqrt{\mathfrak{E}}\int_0^t\mathfrak{D}\ \text{d} s+\text{l.o.t.}.
$$
Thus, in order to absorb this term with the linear part we would have to assume a size condition on the initial data. Instead of doing this, we are going to exploit the structure of this integrand. In particular, after appropriately splitting $I_6$, we will find a perfect derivative and, via an integration by parts, this will allow us to improve the estimates and avoid any size restriction.

We compute
\begin{align*}
I_6&=-\delta\int_0^t\int_{\mathbb{S}^1}\partial_x\left[\partial_x^2f\cH f _t+2\partial_xf\Lambda f _t\right]\Lambda^7 f_t \ \text{d} x \ \text{d}s\\
&=\delta\int_0^t\int_{\mathbb{S}^1}\partial_x^3\left[\partial_x^2f\cH f _t+2\partial_xf\Lambda f _t\right]\Lambda^5 f_t \ \text{d} x \ \text{d}s\\
&=J_1^6+J_2^6+J_3^6+J_4^6+J_5^6,
\end{align*}
with
\begin{align*}
J_1^6&=\delta\int_0^t\int_{\mathbb{S}^1}\partial_x^5f\cH f_t\Lambda^5 f_t \ \text{d} x \ \text{d}s\\
J_2^6&=7\delta\int_0^t\int_{\mathbb{S}^1}\partial_x^2f\Lambda\partial_x^2 f_t\Lambda^5 f_t \ \text{d} x \ \text{d}s\\
J_3^6&=9\delta\int_0^t\int_{\mathbb{S}^1}\partial_x^3f\Lambda\partial_x f_t\Lambda^5 f_t \ \text{d} x \ \text{d}s\\
J_4^6&=5\delta\int_0^t\int_{\mathbb{S}^1}\partial_x^4f\Lambda f_t\Lambda^5 f_t \ \text{d} x \ \text{d}s\\
\end{align*}
and
$$
J_5^6=2\delta\int_0^t\int_{\mathbb{S}^1}\partial_xf\Lambda\partial_x^3 f_t\Lambda^5 f_t \ \text{d} x \ \text{d}s.
$$

To bound the term $J_1^6$ we use the $H^{0.5}-H^{-0.5}$ duality together with the Sobolev embeddging and inequality \eqref{almostalgebra}:
\begin{align*}
J_1^6&\leq \delta\int_0^t\|f\|_{\dot{H}^{5.5}}\|f_t\|_{\dot{H}^1}\|f_t\|_{\dot{H}^{4.5}}\ \text{d} s\\
&\leq \frac{C}{\delta}t\mathfrak{E}^2 \pare{t}+\frac{\delta}{20}\int_0^t\mathfrak{D} \pare{s} \ \text{d} s.
\end{align*}
Similarly,
\begin{align*}
J_2^6&\leq 7\delta\int_0^t\|f\|_{\dot{H}^{3}}\|f_t\|_{\dot{H}^{3.5}}\|f_t\|_{\dot{H}^{4.5}}\ \text{d} s\\
&\leq \frac{C}{\delta}t\mathfrak{E}^2 \pare{t}+\frac{\delta}{20}\int_0^t\mathfrak{D} \pare{s} \ \text{d} s,
\end{align*}
\begin{align*}
J_3^6&\leq 9\delta\int_0^t\|f\|_{\dot{H}^{3.5}}\|f_t\|_{\dot{H}^3}\|f_t\|_{\dot{H}^{4.5}}\ \text{d} s\\
&\leq \frac{C}{\delta}t\mathfrak{E}^2 \pare{t}+\frac{\delta}{20}\int_0^t\mathfrak{D} \pare{s} \ \text{d} s,
\end{align*}
\begin{align*}
J_4^6&\leq 5\delta\int_0^t\|f\|_{\dot{H}^{4.5}}\|f_t\|_{\dot{H}^2}\|f_t\|_{\dot{H}^{4.5}}\ \text{d} s\\
&\leq \frac{C}{\delta}t\mathfrak{E}^2 \pare{t}+\frac{\delta}{20}\int_0^t\mathfrak{D} \pare{s} \ \text{d} s.
\end{align*}
Now we use that
$$
\Lambda^4 f_t=\partial_x^4f_t
$$
to write
$$
J_5^6=2\delta\int_0^t\int_{\mathbb{S}^1}\partial_xf\Lambda\partial_x^3 f_t\Lambda \partial_x^4 f_t \ \text{d} x \ \text{d}s.
$$
Integrating by parts, we find
\begin{align*}
J_5^6&=-\delta\int_0^t\int_{\mathbb{S}^1}\partial_x^2f(\Lambda\partial_x^3 f_t)^2 \ \text{d} x \ \text{d}s\\
&\leq \delta\int_0^t\|f\|_{\dot{H}^{3}}\|f_t\|_{\dot{H}^{4}}^2\ \text{d} s\\
&\leq \delta\int_0^t\|f\|_{\dot{H}^{3}}\|f_t\|_{\dot{H}^{3.5}}\|f_t\|_{\dot{H}^{4.5}}\ \text{d} s\\
&\leq Ct\mathfrak{E}^2\pare{t}+\frac{\delta}{20}\int_0^t\mathfrak{D}\pare{s} \ \text{d} s.
\end{align*}
Thus,
\begin{equation}\label{I62}
I_6\leq  \frac{C}{\delta}t\mathfrak{E}^2\pare{t}+\frac{\delta}{4}\int_0^t\mathfrak{D}\pare{s}\ \text{d} s.
\end{equation}
Collecting the previous estimates for $I^j$ \eqref{I12},\eqref{I22},\eqref{I32},\eqref{I42} and \eqref{I62}, we find
$$
\mathfrak{E}(t)\leq \mathfrak{E}(0)+c(\delta)t(\mathfrak{E}(t))^2.
$$
This last inequality implies the existence of a uniform time such that
$$
\mathfrak{E}(t)\leq 4\mathfrak{E}(0).
$$
Once the estimates are uniform in time, we can mollify \eqref{NEW} as in \cite{granero2017model} to perform a Picard iteration scheme and conclude the existence of a solution
$$
f\in L^\infty\left(0,T;H^{5.5}\right),
$$
$$
f_t\in L^\infty\left(0,T;H^{3.5}\right)\cap L^2(0,T;H^{4.5}).
$$
To obtain the endpoint continuity in time, we observe that one can prove that the solution is continuous in time with respect to the weak topologies. Then, we can perform similar estimates and we obtain a differential inequality that implies the desired continuity in time (see \cite{granero2017model} for more details). Uniqueness of solutions follows from a standard contradiction argument using the $L^2$-type energy estimate, and we omit the details.This finishes the proof.

\begin{remark}
We observe that, to reduce the regularity of the initial data from $H^{5.5}$ to $H^5$, we would have to estimate the term
$$
I=\int_0^t\int_{\mathbb{S}^1}\partial_xf\Lambda \partial_x f _t\Lambda^6 f_t \ \text{d} x \ \text{d}s.
$$
Here, the previous strategy of looking for a perfect derivative in order to integrate by parts will not work. Instead, we would have to use a hidden commutator. Indeed, using that
$$
\int_{\mathbb{S}^1} uv\mathcal{H}v \ \text{d} x =\frac{1}{2}\int_{\mathbb{S}^1} uv\mathcal{H}v \ \text{d} x -\frac{1}{2}\int_{\mathbb{S}^1} \mathcal{H}(uv)v \ \text{d} x =-\frac{1}{2}\int_{\mathbb{S}^1}\comm{ \mathcal{H}}{u}vv \ \text{d} x ,
$$
the term
$$
J=\int_0^t\int_{\mathbb{S}^1}\partial_xf\mathcal{H} \partial_x^4 f _t\partial_x^4 f_t \ \text{d} x \ \text{d}s
$$
has a commutator structure that we can exploit.
\end{remark}

\section*{Acknowledgments}
The research of S.S. is supported by the Basque Government through the BERC 2018-2021 program and by Spanish Ministry of Economy and Competitiveness MINECO through BCAM Severo Ochoa excellence accreditation SEV-2017-0718 and through project MTM2017-82184-R funded by (AEI/FEDER, UE) and acronym "DESFLU" and by the European Research Council through the Starting Grant project H2020-EU.1.1.-639227 FLUID-INTERFACE. R. G-B has been funded by the grant MTM2017-89976-P from the Spanish government.
 
\bibliographystyle{amsplain}

\begin{thebibliography}{10}

\bibitem{AkMi2010}
Benjamin Akers and Paul~A Milewski.
\newblock Dynamics of three-dimensional gravity-capillary solitary waves in
  deep water.
\newblock {\em SIAM Journal on Applied Mathematics}, 70(7):2390--2408, 2010.

\bibitem{AkNi2010}
Benjamin Akers and David~P Nicholls.
\newblock Traveling waves in deep water with gravity and surface tension.
\newblock {\em SIAM Journal on Applied Mathematics}, 70(7):2373--2389, 2010.

\bibitem{ambrose2012well}
David Ambrose, Jerry Bona, and David Nicholls.
\newblock Well-posedness of a model for water waves with viscosity.
\newblock {\em Discrete and Continuous Dynamical Systems, Series B}, 17(4):1113--1137, 2012.

\bibitem{aurther2019rigorous}
CH~Arthur Cheng, Rafael Granero-Belinch{\'o}n, Steve Shkoller, and Jon Wilkening.
\newblock Rigorous asymptotic models of water waves.
\newblock {\em Water Waves}, 1(1):71--130, 2019.

\bibitem{boussinesq1895lois}
J~Boussinesq.
\newblock Lois de l’extinction de la houle en haute mer.
\newblock {\em CR Acad. Sci. Paris}, 121(15-20):2, 1895.

\bibitem{cheng2016well}
CH~Arthur Cheng, Rafael Granero-Belinch{\'o}n, and Steve Shkoller.
\newblock Well-posedness of the Muskat problem with $H^2$ initial data.
\newblock {\em Advances in Mathematics}, 286:32--104, 2016.

\bibitem{craig1993numerical}
Walter Craig and Catherine Sulem.
\newblock Numerical simulation of gravity waves.
\newblock {\em Journal of Computational Physics}, 108(1):73--83, 1993.

\bibitem{dawson2008decay}
L~Dawson, H~McGahagan, and G~Ponce.
\newblock On the decay properties of solutions to a class of Schr{\"o}dinger
  equations.
\newblock {\em Proceedings of the American Mathematical Society},
  136(6):2081--2090, 2008.

\bibitem{dias2008theory}
Frederic Dias, Alexander~I Dyachenko, and Vladimir~E Zakharov.
\newblock Theory of weakly damped free-surface flows: a new formulation based
  on potential flow solutions.
\newblock {\em Physics Letters A}, 372(8):1297--1302, 2008.

\bibitem{dutykh2009visco}
Denys Dutykh.
\newblock Visco-potential free-surface flows and long wave modelling.
\newblock {\em European Journal of Mechanics-B/Fluids}, 28(3):430--443, 2009.


\bibitem{dutykh2007dissipative}
Denys Dutykh and Fr{\'e}d{\'e}ric Dias.
\newblock Dissipative Boussinesq equations.
\newblock {\em Comptes Rendus Mecanique}, 335(9-10):559--583, 2007.

\bibitem{dutykh2007viscous}
Denys Dutykh and Fr{\'e}d{\'e}ric Dias.
\newblock Viscous potential free-surface flows in a fluid layer of finite
  depth.
\newblock {\em Comptes Rendus Mathematique}, 345(2):113--118, 2007.

\bibitem{dutykhA}
Denys Dutykh and Olivier Goubet.
\newblock Derivation of dissipative Boussinesq equations using the
  Dirichlet-to-Neumann operator approach.
\newblock {\em Mathematics and Computers in Simulation}, 127:80--93, 2016.

\bibitem{grafakos2014kato}
Loukas Grafakos and Seungly Oh.
\newblock The Kato-Ponce inequality.
\newblock {\em Communications in Partial Differential Equations},
  39(6):1128--1157, 2014.

\bibitem{granero2018asymptotic}
Rafael Granero-Belinch{\'o}n and Stefano Scrobogna.
\newblock On an asymptotic model for free boundary darcy flow in porous media.
\newblock {\em arXiv preprint arXiv:1810.11798}, 2018.

\bibitem{granero2019asymptotic}
Rafael Granero-Belinch{\'o}n and Stefano Scrobogna.
\newblock Asymptotic models for free boundary flow in porous media.
\newblock {\em Physica D: Nonlinear Phenomena}, 392:1--16, 2019.

\bibitem{graneroscrobo}
Rafael Granero-Belinch{\'o}n and Stefano Scrobogna.
\newblock Models for damped water waves.
\newblock {\em SIAM Journal of Applied Mathematics}, 79:6:2530--2550, 2019.

\bibitem{graneroscrobo2}
Rafael Granero-Belinch{\'o}n and Stefano Scrobogna.
\newblock Well-posedness of the water-wave with viscosity problem.
\newblock {\em Submitted arXiv preprint arXiv:2003.11454}, 2020.

\bibitem{granero2017model}
Rafael Granero-Belinch{\'o}n and Steve Shkoller.
\newblock A model for Rayleigh--Taylor mixing and interface turnover.
\newblock {\em Multiscale Modeling \& Simulation}, 15(1):274--308, 2017.

\bibitem{jiang1996moderate}
Lei Jiang, Chao-Lung Ting, Marc Perlin, and William~W Schultz.
\newblock Moderate and steep Faraday waves: instabilities, modulation and
  temporal asymmetries.
\newblock {\em Journal of Fluid Mechanics}, 329:275--307, 1996.

\bibitem{joseph2004dissipation}
Daniel~D Joseph and Jing Wang.
\newblock The dissipation approximation and viscous potential flow.
\newblock {\em Journal of Fluid Mechanics}, 505:365--377, 2004.

\bibitem{kakleas2010numerical}
Maria Kakleas and David~P Nicholls.
\newblock Numerical simulation of a weakly nonlinear model for water waves with
  viscosity.
\newblock {\em Journal of Scientific Computing}, 42(2):274--290, 2010.

\bibitem{kato1988commutator}
Tosio Kato and Gustavo Ponce.
\newblock Commutator estimates and the Euler and Navier-Stokes equations.
\newblock {\em Communications on Pure and Applied Mathematics}, 41(7):891--907,
  1988.

\bibitem{kenig1993well}
Carlos~E Kenig, Gustavo Ponce, and Luis Vega.
\newblock Well-posedness and scattering results for the generalized Korteweg-de
  Vries equation via the contraction principle.
\newblock {\em Communications on Pure and Applied Mathematics}, 46(4):527--620,
  1993.

\bibitem{khariff1996frequency}
C~Khariff, Skandrani C, and J~Poitevin.
\newblock The frequency down-shift phenomenon.
\newblock In {\em Mathematical Problems in the Theory of Water Waves: A
  Workshop on the Problems in the Theory of Nonlinear Hydrodynamic Waves, May
  15-19, 1995, Luminy, France}, volume 200, page 157. American Mathematical
  Soc., 1996.

\bibitem{lamb1932hydrodynamics}
H~Lamb.
\newblock {\em Hydrodynamics}.
\newblock Cambridge Univ Press,, 1932.

\bibitem{longuet1992theory}
Michael~S Longuet-Higgins.
\newblock Theory of weakly damped Stokes waves: a new formulation and its
  physical interpretation.
\newblock {\em Journal of Fluid Mechanics}, 235:319--324, 1992.

\bibitem{matsuno1992nonlinear}
Y~Matsuno.
\newblock Nonlinear evolutions of surface gravity waves on fluid of finite
  depth.
\newblock {\em Physical review letters}, 69(4):609, 1992.

\bibitem{matsuno1993nonlinear}
Yoshimasa Matsuno.
\newblock Nonlinear evolution of surface gravity waves over an uneven bottom.
\newblock {\em Journal of fluid mechanics}, 249:121--133, 1993.

\bibitem{matsuno1993two}
Yoshimasa Matsuno.
\newblock Two-dimensional evolution of surface gravity waves on a fluid of
  arbitrary depth.
\newblock {\em Physical Review E}, 47(6):4593, 1993.

\bibitem{ngom2018well}
Mari{\`e}me Ngom and David~P Nicholls.
\newblock Well-posedness and analyticity of solutions to a water wave problem
  with viscosity.
\newblock {\em Journal of Differential Equations}, 265(10):5031--5065, 2018.

\bibitem{ruvinsky1991numerical}
KD~Ruvinsky, FI~Feldstein, and GI~Freidman.
\newblock Numerical simulations of the quasi-stationary stage of ripple
  excitation by steep gravity--capillary waves.
\newblock {\em Journal of Fluid Mechanics}, 230:339--353, 1991.

\bibitem{ruvinsky1985improvement}
KD~Ruvinsky and GI~Freidman.
\newblock Improvement of the first Stokes method for the investigation of
  finite-amplitude potential gravity-capillary waves.
\newblock In {\em IX All-Union Symp. on Diffraction and Propagation Waves,
  Tbilisi: Theses of Reports}, volume~2, pages 22--25, 1985.

\bibitem{ruvinsky1987fine}
KD~Ruvinsky and GI~Freidman.
\newblock The fine structure of strong gravity-capillary waves.
\newblock {\em Nonlinear waves: Structures and Bifurcations, AV Gaponov-Grekhov
  and MI Rabinovich, eds. Moscow: Nauka}, pages 304--326, 1987.

\bibitem{Stokes_1847}
George~Gabriel Stokes.
\newblock {\em On the Theory of Oscillatory Waves}, volume~1 of {\em Cambridge
  Library Collection - Mathematics}.
\newblock Cambridge University Press, 1847.

\bibitem{wang2006purely}
Jing Wang and Daniel~D Joseph.
\newblock Purely irrotational theories of the effect of the viscosity on the
  decay of free gravity waves.
\newblock {\em Journal of Fluid Mechanics}, 559:461--472, 2006.

\bibitem{wu2006note}
Guangyu Wu, Yuming Liu, and Dick~KP Yue.
\newblock A note on stabilizing the Benjamin--Feir instability.
\newblock {\em Journal of Fluid Mechanics}, 556:45--54, 2006.

\end{thebibliography}

\end{document}